\documentclass[a4paper,oneside,12pt]{article}
\usepackage{amssymb}

\usepackage{amsmath}
\usepackage{a4}

\input{tcilatex}

\begin{document}

\title{A unified formula for Steenrod operations in flag manifolds}
\author{Haibao Duan and Xuezhi Zhao \\
Institute of Mathematics, Chinese Academy of Sciences, \\
Beijing 100080, dhb@math.ac.cn\\
Department of Mathematics, Capital Normal University\\
Beijing 100037, zhaoxve@mail.cnu.edu.cn}
\date{\ \ }
\maketitle

\begin{abstract}
The classical Schubert cells on a flag manifold $G/H$ give a cell
decomposition for $G/H$ whose Kronecker duals (known as Schubert classes)
form an additive base for the integral cohomology $H^{\ast }(G/H)$.

We present a formula that expresses Steenrod mod--p operations on Schubert
classes in $G/H$ in terms of Cartan numbers of $G$.

\begin{description}
\item[2000 Mathematical Subject Classification:] 55S10 (14M10, 1415).

\item[Key words and Phreses:] Steenrod operations, Flag manifolds, Schubert
cells, Cartan numbers.
\end{description}
\end{abstract}

\begin{center}
\textbf{1. Introduction}
\end{center}

Let $p\geq 2$ be a fixed prime and let $\mathcal{A}_{p}$ be the mod--$p$
Steenrod algebra. Denote by $\mathcal{P}^{k}\in \mathcal{A}_{p}$, $k\geq 0$,
the Steenrod mod--$p$ reduced powers on the $\mathbb{Z}_{p}$--cohomology of
topological spaces [SE]. When $p=2$, it is more customary to write $%
Sq^{2k}\in \mathcal{A}_{2}$ instead of $\mathcal{P}^{k}$.

In general, an entire description of the $\mathcal{A}_{p}$--action on the $\mathbb{Z}%
_{p} $--cohomology of a topological space $X$ leads to two enquiries.

\begin{quote}
\textbf{Problem A.}\textsl{\ Specify an additive basis }$S=\{\omega
_{1},\cdots ,\omega _{m}\}$\textsl{\ for the graded }$\mathbb{Z}_{p}$\textsl{%
--vector space }$H^{\ast }(X;\mathbb{Z}_{p})$\textsl{\ that encodes the
geometric formation of }$X$\textsl{\ (e.g. a cell decomposition of }$X$%
\textsl{).}

\textbf{Problem B.} \textsl{Determine all the coefficients }$c_{k,i}^{j}\in
Z_{p}$\textsl{\ in the expression}
\end{quote}

\begin{center}
$P^{k}(\omega _{i})=\sum_{1\leq j\leq m}c_{k,i}^{j}\omega _{j}$\textsl{, }$%
k\geq 0$\textsl{, }$1\leq i\leq m$\textsl{.}
\end{center}

The study of the internal structure of the algebra $\mathcal{A}_{p}$ has
continued for almost fifty years, see Wood [W] for a thorough historical
account, further problems and relevant references. On the other hand, we
note that even partial solutions to Problem B can have significant
consequences in manifold geometry. To mention only a few examples: the
classical Wu--formula [Wu] can be interpreted as the expansion of the $%
Sq^{k} $--action on the special Schubert classes in the real Grassmannians;
the calculation by Steenrod and Whitehead [SW] in the truncated real
projective spaces led to an enormous step in understanding the classical
problem of how many linearly independent vector fields can be found on the $%
n $-sphere $S^{n} $; by deriving partial knowledge of the $\mathcal{P}^{k}$%
--action on the special Schubert classes in the complex Grassmannian, Borel
and Serre ([BSe$_{1}$,BSe$_{2}$]) demonstrated that the $2n$--dimensional
spheres $S^{2n}$ does not admit any almost complex structure\ unless $%
n=1,2,3 $. Needless to say many profound applications and deep implication
of Steenrod operations in topology [D,St,Le], effective computation of these
operations in given manifolds deserves also high priority.

\bigskip\

Let $G$ be a compact connected Lie group and let $H$ be the centralizer of a
one-parameter subgroup in $G$. The space $G/H=\{gH\mid g\in G\}$ of left
cosets of $H$ in $G$ is known as a \textsl{flag manifold}. In this paper we
study Problem A and B for all $G/H$ .

Firstly, if $X$ is a flag manifold $G/H$, a uniform solution to Problem A is
already known as\textbf{\ }\textsl{Basis Theorem}\textbf{\ }from \textsl{the
Schubert's enumerative calculus} (i.e. a branch of algebraic geometry [K,
So, BGG]). It was originated by Ehresmann [E] for the Grassmannians $G_{n,k}$
of $k$-dimensional subspaces in $\mathbb{C}^{n}$ in 1934, extended to the
case where $G$ is a matrix group by Bruhat in 1954, and completed for all
compact connected Lie groups by Chevalley [Ch] in 1958. We briefly recall
the result.

Let $W$ and $W^{^{\prime }}$ be the Weyl groups of $G$ and $H$ respectively.
The set $W/W^{^{\prime }}$ of left cosets of $W^{^{\prime }}$ in $W$ can be
identified with the subset of $W$:

\begin{center}
$\overline{W}=\{w\in W\mid l(w_{1})\geq l(w)$ for all $w_{1}\in wW^{^{\prime
}}\}$,\
\end{center}

\noindent where $l:W\rightarrow\mathbb{Z}$ is the length function relative
to a fixed maximal torus $T$ in $G$ [BGG, 5.1. Proposition]. The key fact is
that the space $G/H$ admits a canonical decomposition into cells indexed by
elements of $\overline{W}$

\begin{enumerate}
\item[(1.1)] $\ \ \ \ \ \ \ \ G/H=\underset{w\in \overline{W}}{\cup }X_{w}$,
$\quad \dim X_{w}=2l(w)$,
\end{enumerate}

\noindent with each cell $X_{w}$ the closure of an algebraic affine space,
known as a \textsl{Schubert variety} in $G/H$ [Ch, BGG]. Since only even
dimensional cells are involved in the decomposition (1.1), the set of
fundamental classes\textsl{\ }$[X_{w}]\in H_{2l(w)}(G/H)$, $w\in \overline{W}
$,\textsl{\ }forms an additive basis of the homology $H_{\ast }(G/H)$. The
cocycle class $P_{w}(H)\in H^{2l(w)}(G/H)$ defined by the Kronecker pairing
as

\begin{center}
$\left\langle P_{w}(H),[X_{u}]\right\rangle =\delta _{w,u}$, $w,u\in
\overline{W}$,
\end{center}

\noindent is called the \textsl{Schubert class corresponding to }$w$.
Combining (1.1) with the Poincar\'{e} duality yields the following solution
to Problem A.

\textbf{Lemma 1 }(Basis Theorem\textbf{).} \textsl{The set of Schubert
classes }$\{P_{w}(H)\mid $\textsl{\ }$w\in \overline{W}\}$\textsl{\
constitutes an additive basis for the ring }$H^{\ast }(G/H)$\textsl{.}

\bigskip

It follows from the Basis Theorem that, for a $u\in \overline{W}$ and $k\geq
0$, one has the expression

\begin{enumerate}
\item[(1.2)] $\ \ \ \ \ \ \ \ \mathcal{P}^{k}(P_{u}(H))\equiv \sum
a_{w,u}^{k}P_{w}(H)\mod p$, $a_{w,u}^{k}\in \mathbb{Z}_{p}$,
\end{enumerate}

\noindent where the sum ranges over all $w\in \overline{W}$ with $%
l(w)=l(u)+k(p-1)$ for dimension reason. Thus, in the case of $X=G/H$,
Problem B admits a concrete form.

\begin{quote}
\textbf{Problem C.} \textsl{Determine the numbers }$a_{w,u}^{k}\in \mathbb{Z}%
_{p}$\textsl{\ for }$k\geq 0$\textsl{, }$w,u\in \overline{W}$ \textsl{with }$%
l(w)=l(u)+k(p-1)$\textsl{.}
\end{quote}

If $G$ is the unitary group $U(n)$ of order $n$ and $H=U(k)\times U(n-k$),
the flag manifold $G/H$ is the complex Grassmannian $G_{n,k}$ of $k$-planes
through the origin in $\mathbb{C}^{n}$. The $i^{th}$ Chern classes $c_{i}\in
H^{2i}(G_{n,k})$, $1\leq i\leq k$, of the canonical complex $k$-bundle over $%
G_{n,k}$ are precisely \textsl{the special Schubert classes} on $G_{n,k}$
[So, GH]. In order to generalize the classical Wu--formula [Wu], many works
were devoted to find an expression of $\mathcal{P}^{k}(c_{i})$ in terms of
the $c_{i}$ (cf. [BSe$_{1}$, BSe$_{2}$, D$_{4}$, La, Le, P, S, Su]). This
seems to be the only special case for which Problem C has been studied in
some details.

It is well known that the knowledge of the $\mathcal{A}_{p}$--action on the $%
\mathbb{Z}_{p}$--cohomology of a space $X$ can provide deeper information on
the topology of $X$ than just the cohomology ring structure. In this regard
the present work is as a sequel of [D$_{3}$,DZ], where the multiplicative
rule of Schubert classes (in the integral cohomology $H^{\ast }(G/H)$) was
determined. In particular, in view of the geometric decomposition (1.1) of $%
G/H$ offered by the classical Schubert cells, the numbers $a_{w,u}^{k}$ are
immediately applicable to investigate the attaching maps of these cells
(e.g. Compare the tables in Section 5 with the figure in [Le, Section 6]).
We quote from Lenart [Le] for the case of $G_{n,k}$: \textsl{Apart from
projective spaces, very little is known about the attaching maps of their
cells}.

\bigskip

This paper is arranged as follows. Section 2 contains a brief introduction
to the Weyl group associated with a Lie group. Then the solution to Problem
C (i.e. the Theorem) is presented. After geometric preliminaries in Section
3 the Theorem is established in Section 4. In order to illustrate the
effective computability of our method, computational results for some cases
of $G/H$ are explained and tabulated in the final section.

\begin{center}
\textbf{2. The result}
\end{center}

To investigate a flag manifold $G/H$ one may assume that the Lie group $G$
under consideration is $1$--connected and semi-simple ([BH]). Since all $1$%
--connected semi--simple Lie groups are classified by their Cartan matrices
[Hu, p.55], any numerical topological invariant of $G/H$ \ may be reduced to
the Cartan numbers (entries in the Cartan matrix of $G$). We present both a
formula and an algorithm, that evaluate the numbers $a_{w,u}^{k}\in \mathbb{Z%
}_{p}$ in terms of Cartan numbers of $G$.

\bigskip

Fix a maximal torus $T$ of $G$ and set $n=\dim T$. Equip the Lie algebra $%
L(G)$ of $G$ with an inner product $(,)$ so that the adjoint representation
acts as isometries of $L(G)$. The\textsl{\ Cartan subalgebra} of $G$ is the
Euclidean subspace $L(T)$ of $L(G)$.

The restriction of the exponential map $\exp :L(G)\rightarrow G$ to $L(T)$
defines a set $D(G)$ of $\frac{1}{2}(\dim G-n)$ hyperplanes in $L(T)$, i.e.
the set of\textsl{\ singular hyperplanes }through the origin in $L(T)$. The
reflections $\sigma $ of $L(T)$ in these planes generate the Weyl group $W$
of $G$ ([Hu, p.49]).

Fix a regular point $\alpha \in L(T)\backslash \cup _{L\in D(G)}L$ and let $%
\Delta $ be the set of simple roots relative to $\alpha $ [Hu, p.47]. If $%
\beta \in \Delta $ the reflection $\sigma _{\beta }$ in the hyperplane $%
L_{\beta }\in D(G)$ relative to $\beta $ is called a \textsl{simple
reflection}. If $\beta ,\beta ^{\prime }\in \Delta $, \textsl{the Cartan
number }

\begin{center}
$\beta \circ \beta ^{\prime }=2(\beta ,\beta ^{\prime })/(\beta ^{\prime
},\beta ^{\prime })$
\end{center}

\noindent is always an integer (only $0,\pm 1,\pm 2,\pm 3$ can occur [Hu,
p.55]).

\bigskip

It is known that the set of simple reflections $\{\sigma _{\beta }\mid \beta
\in \Delta \}$ generates $W$. That is, any $w\in W$ admits a factorization
of the form

\begin{enumerate}
\item[(2.1)] \ \ \ \ \ \ \noindent $w=\sigma _{\beta _{1}}\circ \cdots \circ
\sigma _{\beta _{m}}$, $\beta _{i}\in \Delta $.
\end{enumerate}

\textbf{Definition 1}. The \textsl{length} $l(w)$ of a $w\in W$ is the least
number of factors in all decompositions of $w$ in the form (2.1). The
decomposition (2.1) is said \textsl{reduced} if $m=l(w)$.

If (2.1) is a reduced decomposition, the $m\times m$ (strictly upper
triangular) matrix $A_{w}=(a_{i,j})$ with

\begin{center}
$a_{i,j}=\{%
\begin{array}{c}
0\text{ if }i\geq j\text{;\qquad } \\
-\beta _{i}\circ \beta _{j}\text{ if }i<j%
\end{array}%
$
\end{center}

\noindent will be called \textsl{the Cartan matrix of }$w$\textsl{\ }%
associated to the decomposition\textsl{\ }(2.1).

\bigskip

Let $\mathbb{Z}[x_{1},\cdots ,x_{m}]=\oplus _{n\geq 0}\mathbb{Z}%
[x_{1},\cdots ,x_{m}]^{(n)}$ be the ring of integral polynomials in $%
x_{1},\cdots ,x_{m}$, graded by $\mid x_{i}\mid =1$.

\textbf{Definition 2.} For a subset $[i_{1},\cdots ,i_{r}]\subseteq \lbrack
1,\cdots ,m]$ and $1\leq k\leq r$, denote by $m_{k,p}(x_{i_{1}},\cdots
,x_{i_{r}})$ the polynomial

\begin{center}
$\underset{(\alpha _{1},\cdots ,\alpha _{r})}{\Sigma }x_{i_{1}}^{\alpha
_{1}}\cdots x_{i_{r}}^{\alpha _{r}}\in \mathbb{Z}[x_{1},\cdots
,x_{m}]^{(r+k(p-1))}$,
\end{center}

\noindent where the sum is over all distinct permutations $(\alpha
_{1},\cdots ,\alpha _{r})$ of the partition $(p^{k},1^{r-k})$ ([M, p.1]).

\bigskip

\textbf{Remark 1.} In the theory of symmetric functions, the $%
m_{k,p}(x_{i_{1}},\cdots ,x_{i_{r}})$ is known as \textsl{the monomial
symmetric function }in $x_{i_{1}},\cdots ,x_{i_{r}}$ associated to the
partition $(p^{k},1^{r-k})$ ([M, p.19]). As examples, if $[i_{1},\cdots
,i_{r}]=[1,2,3]$ one has

\begin{quote}
$m_{1,p}(x_{i_{1}},\cdots
,x_{i_{r}})=x_{1}^{p}x_{2}x_{3}+x_{1}x_{2}^{p}x_{3}+x_{1}x_{2}x_{3}^{p}$;

$m_{2,p}(x_{i_{1}},\cdots
,x_{i_{r}})=x_{1}^{p}x_{2}^{p}x_{3}+x_{1}x_{2}^{p}x_{3}^{p}+x_{1}^{p}x_{2}x_{3}^{p}
$;

$m_{3,p}(x_{i_{1}},\cdots ,x_{i_{r}})=x_{1}^{p}x_{2}^{p}x_{3}^{p}$.
\end{quote}

\bigskip

\textbf{Definition 3. }Given a $m\times m$ strictly upper triangular integer
matrix $A=(a_{i,j})$ define a homomorphism $T_{A}:$ $\mathbb{Z}[x_{1},\cdots
,x_{m}]^{(m)}\rightarrow \mathbb{Z}$ recursively as follows:

\begin{quote}
1) for $h\in \mathbb{Z}[x_{1},\cdot \cdot \cdot ,x_{m-1}]^{(m)}$, $%
T_{A}(h)=0 $;

2) if $m=1$ (consequently $A=(0)$), then $T_{A}(x_{1})=1$;

3) for $h\in \mathbb{Z}[x_{1},\cdot \cdot \cdot ,x_{m-1}]^{(m-r)}$ with $%
r\geq 1$,
\end{quote}

\begin{center}
$T_{A}(hx_{m}^{r})=T_{A^{\prime }}(h(a_{1,m}x_{1}+\cdots
+a_{m-1,m}x_{m-1})^{r-1})$,
\end{center}

\begin{quote}
\noindent where $A^{\prime }$ is the ($(m-1)\times (m-1)$ strictly upper
triangular) matrix obtained from $A$ by deleting the $m^{th}$ column and the
$m^{th}$ row.
\end{quote}

\noindent By additivity, $T_{A}$ is defined for every $h\in \mathbb{Z}%
[x_{1},\cdots ,x_{m}]^{(m)}$ using the unique expansion $h=\sum\limits_{0%
\leq r\leq m}h_{r}x_{m}^{r}$ with $h_{r}\in \mathbb{Z}[x_{1},\cdot \cdot
\cdot ,x_{m-1}]^{(m-r)}$.

\bigskip

\textbf{Remark 2.} Definition 3 implies an effective algorithm to evaluate $%
T_{A}$.

For $k=2$ and $A_{1}=\left(
\begin{array}{cc}
0 & a \\
0 & 0%
\end{array}%
\right) $, then $T_{A_{1}}:$ $\mathbb{Z}[x_{1},x_{2}]^{(2)}\rightarrow
\mathbb{Z}$ is given by

\begin{quote}
$T_{A_{1}}(x_{1}^{2})=0$,$\qquad $

$T_{A_{1}}(x_{1}x_{2})=T_{A_{1}^{\prime }}(x_{1})=1$ and

$T_{A_{1}}(x_{2}^{2})=T_{A_{1}^{\prime }}(ax_{1})=a$.
\end{quote}

\noindent For $k=3$ and $A_{2}=\left(
\begin{array}{ccc}
0 & a & b \\
0 & 0 & c \\
0 & 0 & 0%
\end{array}%
\right) $, then $A_{2}^{\prime }=A_{1}$ and $T_{A_{2}}:$ $\mathbb{Z}%
[x_{1},x_{2},x_{3}]^{(3)}\rightarrow \mathbb{Z}$ is given by

\begin{center}
$T_{A_{2}}(x_{1}^{r_{1}}x_{2}^{r_{2}}x_{3}^{r_{3}})=\{%
\begin{array}{c}
0\text{, if }r_{3}=0\text{ and\qquad \qquad \qquad \qquad \qquad } \\
T_{A_{1}}(x_{1}^{r_{1}}x_{2}^{r_{2}}(bx_{1}+cx_{2})^{r_{3}-1}),\text{ if }%
r_{3}\geq 1\text{,\ }%
\end{array}%
$
\end{center}

\noindent where $r_{1}+r_{2}+r_{3}=3$, and where $T_{A_{1}}$ is calculated
in the above.

\bigskip

\bigskip

Assume that $w=\sigma _{\beta _{1}}\circ \cdots \circ \sigma _{\beta _{m}}$,
$\beta _{i}\in \Delta $, is a reduced decomposition of $w\in \overline{W}$
and let $A_{w}=(a_{i,j})_{m\times m}$ be the associated Cartan matrix of $w$%
. For a subset $J=[i_{1},\cdots ,i_{r}]\subseteq \lbrack 1,\cdots ,m]$ we set

\begin{center}
$\sigma _{J}=\sigma _{\beta _{i_{1}}}\circ \cdots \circ \sigma _{\beta
_{i_{r}}}$.
\end{center}

\noindent Our solution to Problem C is

\textbf{Theorem.} \textsl{For a }$u\in \overline{W}$\textsl{,} $k>0$ \textsl{%
with }$l(u)=r$\textsl{\ and }$m=r+k(p-1)$\textsl{, we have (in (1.2)) that}

\begin{center}
$\qquad a_{w,u}^{k}\equiv T_{A_{w}}[\sum\limits_{\substack{ J=[i_{1},\cdots
,i_{r}]\subseteq \lbrack 1,\cdots ,m]  \\ \sigma _{J}=u}}m_{k,p}(x_{i_{1}},%
\cdots ,x_{i_{r}})]$ \textsl{\ mod p}\textsl{.}
\end{center}

This result will be shown in Section 4.

\bigskip

Indeed, the Theorem indicates an effective algorithm to evaluate $%
a_{w,u}^{k} $ as the following recipe shows.

\begin{quote}
(1) Starting from the Cartan matrix of $G$, \ a program to enumerate all
elements in $\overline{W}$ by their \textsl{minimal reduced decompositions}
is available in [DZ], [DZZ];

(2) For a $w\in \overline{W}$ with a reduced decomposition, the
corresponding Cartan matrix $A_{w}$ can be read directly from Cartan matrix
of $G$ (Compare Definition 1 with [Hu, p.59]);

(3) For a $w\in \overline{W}$ with a reduced decomposition $w=\sigma _{\beta
_{1}}\circ \cdots \circ \sigma _{\beta _{m}}$ and a $u\in \overline{W}$ with
$l(u)=r<m$, the solutions in the subsequence $J=[i_{1},\cdots
,i_{r}]\subseteq \lbrack 1,\cdots ,m]$ to the equation $\sigma _{J}=u$ in $W$
agree with the solutions to the equation $\sigma _{J}(\alpha )=u(\alpha )$
in the vector space $L(T)$, where $\alpha \in L(T)$ is a fixed regular point;

(4) The evaluation the operator $T_{A_{w}}$ on a polynomial can be easily
programmed (cf. Definition 3 or [DZ, Section 5]).
\end{quote}

Based on the algorithm explained above, a program to produce the numbers $%
a_{w,u}^{k}$ has been compiled. As examples some computational results from
the program are explained and tabulated in Section 5.

\bigskip

\bigskip

\begin{center}
\textbf{3. Preliminaries in the }$K$\textbf{--cycles of Bott--Samelson}
\end{center}

In this section all homologies (resp. cohomologies) will have integer
coefficients unless otherwise specified. If $f:X\rightarrow Y$ is a
continuous map between two topological spaces, $f_{\ast }$ (resp. $f^{\ast }$%
) is the homology (resp. cohomology) map induced by $f$. If $M$ is an
oriented closed manifold (resp. a connected projective variety) $[M]\in
H_{\dim M}(M)$ stands for the orientation class. The Kronecker pairing,
between cohomology and homology of a space $X$, will be denoted by $%
<,>:H^{\ast }(X)\times H_{\ast }(X)\rightarrow \mathbb{Z}$.

The proof of our Theorem will make use of the celebrated $K$-cycles (i.e.
\textsl{Bott-Samelson resolutions of Schubert varieties}) on the flag
manifold $G/T$ constructed by Bott and Samelson early in 1955 [BS$_{1}$]. In
this Section we recall the construction of these cycles, as well as their
basic properties (from Lemma 2 to Lemma 4) developed in [BS$_{1}$,BS$_{2}$,D$%
_{1}$,D$_{3}$]. The main technique result in this section is Lemma 5, which
allows us to transform the proof of the Theorem for $G/H$ to calculation in
the $K$--cycles of Bott--Samelson.

\bigskip

As in Section 2, we fix a regular point $\alpha \in L(T)$ and let $\Delta $
be the set of simple roots relative to $\alpha $. For a $\beta \in \Delta $,
the singular plane in $L(T)$ relative to $\beta $ will be denoted by $%
L_{\beta }$ [Hu,p.47]. Write by $K_{\beta }$ the centralizer of $\exp
(L_{\beta })$ in $G$, where $\exp $ is the restriction of the exponential
map $L(G)\rightarrow G$ to $L(T)$. We note that $T\subset K_{\beta }$ and
that the quotient manifold $K_{\beta }/T$ is diffeomorphic to $2$-sphere [BS$%
_{2}$, p.996].

The $2$-sphere $K_{\beta }/T$ carries a natural orientation $\omega _{\beta
}\in H_{2}(K_{\beta }/T;\mathbb{Z})$ that may be specified as follows. The
Cartan decomposition of the Lie algebra $L(K_{\beta })$ relative to the
maximal torus $T\subset K_{\beta }$ has the form $L(K_{\beta })=L(T)\oplus
\vartheta _{\beta }$, where $\vartheta _{\beta }\subset L(G)$ is a $2$%
-plane, the \textsl{root space} belonging to the root $\beta $ [Hu, p.35].
Let $[,]$ be the Lie bracket on $L(G)$. Take a non-zero vector $v\in
\vartheta _{\beta }$ and let $v^{\prime }\in \vartheta _{\beta }$ be such
that $[v,v^{\prime }]=\beta $. The ordered base $\{v,v\prime \}$ gives an
orientation on $\vartheta _{\beta }$ which does not depend on the initial
choice of $v$.

The tangential of the quotient map $\pi _{\beta }:K_{\beta }\rightarrow
K_{\beta }/T$ at the group unit $e\in K_{\beta }$ maps $\vartheta _{\beta }$
isomorphically onto the tangent space to $K_{\beta }/T$ at $\pi _{\beta }(e)$%
. In this manner the orientation $\{v,v\prime \}$ on $\vartheta _{\beta }$
furnishes $K_{\beta }/T$ with the induced orientation $\omega _{\beta
}=\{\pi _{\beta }(v),\pi _{\beta }(v\prime )\}$.

\bigskip

For a sequence $\beta _{1},\cdots ,\beta _{m}\in \Delta $ of simple roots
(repetition like $\beta _{i}=\beta _{j}$ may occur) let $K(\beta _{1},\cdots
,\beta _{m})$ be the product group $K_{\beta _{1}}\times \cdots \times
K_{\beta _{m}}$. Since $T\subset K_{\beta _{i}}$ for each $i$ the group $%
T\times \cdots \times T$ ($m$-copies) acts on $K(\beta _{1},\cdots ,\beta
_{m})$ from the right by

\begin{center}
$(g_{1},\cdots ,g_{m})(t_{1},\cdots
,t_{m})=(g_{1}t_{1},t_{1}^{-1}g_{2}t_{2},\cdots ,t_{m-1}^{-1}g_{m}t_{m})$.
\end{center}

\noindent Let $\Gamma (\beta _{1},\cdots ,\beta _{m})$ be the base manifold
of this principle action, oriented by the $\omega _{\beta _{i}}$, $1\leq
i\leq m$. The point in $\Gamma (\beta _{1},\cdots ,\beta _{m})$
corresponding to a $(g_{1},\cdots ,g_{m})\in K(\beta _{1},\cdots ,\beta
_{m}) $ is denoted by $[g_{1},\cdots ,g_{m}]$.

The integral cohomology of $\Gamma (\beta _{1},\cdots ,\beta _{m})$ has been
determined in [BS$_{1}$, Proposition II]. Let $\varphi _{i}:K_{\beta
_{i}}/T\rightarrow \Gamma (\beta _{1},\cdots ,\beta _{m})$ be the embedding
induced by the inclusion $K_{\beta _{i}}\rightarrow K(\beta _{1},\cdots
,\beta _{m})$ onto the $i^{th}$ factor group, and put

\begin{center}
$y_{i}=\varphi _{i\ast }(\omega _{\beta _{i}})\in H_{2}(\Gamma (\beta
_{1},\cdots ,\beta _{m}))$, $1\leq i\leq m$.
\end{center}

\noindent Form the $m\times m$ strictly upper triangular matrix $%
A=(a_{i,j})_{m\times m}$ by letting

\begin{center}
$a_{i,j}=\{%
\begin{array}{c}
0\text{ if }i\geq j\text{;\qquad } \\
-\beta _{i}\circ \beta _{j}\text{ if }i<j\text{.}%
\end{array}%
$
\end{center}

\textbf{Lemma 2 }([BS$_{2}$])\textbf{.} \textsl{The set }$\{y_{1},\cdots
,y_{m}\}$\textsl{\ forms a basis for }$H_{2}(\Gamma (\beta _{1},\cdots
,\beta _{m}))$\textsl{. Further, let }$x_{i}\in H^{2}(\Gamma (\beta
_{1},\cdots ,\beta _{m}))$\textsl{, }$1\leq i\leq m$\textsl{, be the classes
Kronecker dual to }$y_{1},\cdots ,y_{m}$\textsl{\ as }$<x_{i},y_{j}>=\delta
_{i,j}$\textsl{, }$1\leq i,j\leq m$\textsl{,then}

\begin{center}
$H^{\ast }(\Gamma (\beta _{1},\cdots ,\beta _{m}))=\mathbb{Z}[x_{1},\cdots
,x_{m}]/I$,
\end{center}

\noindent \textsl{where }$I$\textsl{\ is the idea generated by }

\begin{center}
$x_{j}^{2}-\underset{i<j}{\Sigma }a_{i,j}x_{i}x_{j}$\textsl{, }$1\leq j\leq
m $\textsl{.}
\end{center}

\bigskip

In view of Lemma 2 we introduce an additive map $\int_{\Gamma (\beta
_{1},\cdots ,\beta _{m})}:\mathbb{Z}[x_{1},\cdot \cdot \cdot
,x_{m}]^{(m)}\rightarrow \mathbb{Z}$ by

\begin{center}
$\int_{\Gamma (\beta _{1},\cdots ,\beta _{m})}h=<p_{\Gamma (\beta
_{1},\cdots ,\beta _{m})}(h),[\Gamma (\beta _{1},\cdots ,\beta _{m})]>$,
\end{center}

\noindent where $[\Gamma (\beta _{1},\cdots ,\beta _{m})]\in H_{2m}(\Gamma
(\beta _{1},\cdots ,\beta _{m}))=\mathbb{Z}$ is the orientation\textsl{\ }%
class and where

\begin{center}
$p_{\Gamma (\beta _{1},\cdots ,\beta _{m})}:\mathbb{Z}[x_{1},\cdots
,x_{m}]\rightarrow H^{\ast }(\Gamma (\beta _{1},\cdots ,\beta _{m}))$
\end{center}

\noindent is the obvious quotient homomorphism. The geometric implication of
the operator $T_{A}$ in Definition 3 (Section 2) is seen from the next
result.

\textbf{Lemma 3 }([D$_{1}$, Proposition 2]). \textsl{We have}

\begin{center}
\textsl{\qquad \qquad\ \ }$\int_{\Gamma (\beta _{1},\cdots ,\beta
_{m})}=T_{A}:\mathbb{Z}[x_{1},\cdots ,x_{m}]^{(m)}\rightarrow \mathbb{Z}$%
\textsl{.\ }
\end{center}

\noindent \textsl{In particular, }$\int_{\Gamma (\beta _{1},\cdots ,\beta
_{m})}x_{1}\cdots x_{m}=1$\textsl{.}

\bigskip

It follows also from Lemma 2 that the ring $H^{\ast }(\Gamma (\beta
_{1},\cdots ,\beta _{m}))$ has the additive basis $\{x_{i_{1}}\cdots
x_{i_{r}}\mid \lbrack i_{1},\cdots ,i_{r}]\subseteq \lbrack 1,\cdots ,m]\}$.
Since the dimension of every $x_{i}$ is $2$, the action of the $\mathcal{P}%
^{k}\in \mathcal{A}_{p}$ on these base elements is determined by the Cartan
formula [SE]. Let $m_{k,p}(x_{i_{1}},\cdots ,x_{i_{r}})$\ be the monomial
symmetric function\ in $x_{i_{1}},\cdots ,x_{i_{r}}$\ associated to the
partition $(p^{k},1^{r-k})$\textsl{\ }(cf. Definition 2).

\textbf{Lemma 4.} $\mathcal{P}^{k}(x_{i_{1}}\cdots x_{i_{r}})\equiv
m_{k,p}(x_{i_{1}},\cdots ,x_{i_{r}})$\textsl{\ mod p}.

\bigskip

Let $H$ be the centralizer of a one-parameter subgroup in $G$ and let $G/H$
be the flag manifold of left cosets of $H$ in $G$. Assume (without loss the
generality) that with respect to the fixed maximal torus $T\subset G$, $%
T\subseteq H\subset G$.

\textbf{Definition 4.} The map

\begin{center}
$\varphi _{\beta _{1},\cdots ,\beta _{m};H}:\Gamma (\beta _{1},\cdots ,\beta
_{m})\rightarrow G/H$
\end{center}

\noindent by $[g_{1},\cdots ,g_{m}]\rightarrow g_{1}\cdots g_{m}H$ is
clearly well defined and will be called the \textsl{K-cycle} of
Bott-Samelson on $G/H$ \ associated to the sequence $\beta _{1},\cdots
,\beta _{m}\in \Delta $ of simple roots (cf. [D$_{3}$, Subsection 7.1]).

\bigskip

It was first shown by Hansen [H] in 1972 that, when $H=T$, certain K-cycles
of Bott--Samelson provide disingulations of Schubert varieties in $G/T$. The
following more general result allows us to bring the calculation of $%
\mathcal{P}^{k}$--action on the ring $H^{\ast }(G/H)$ (i.e. Problem C) to
the computation of the action on the truncated polynomial algebra $H^{\ast
}(\Gamma (\beta _{1},\cdots ,\beta _{m}))$, while the latter is easily
handled by Lemma 4.

\bigskip

\textbf{Lemma 5.} \textsl{Let }$\{P_{w}(H)\in H^{2l(w)}(G/H)\mid $\textsl{\ }%
$w\in \overline{W}\}$\textsl{\ be the set of Schubert classes on }$G/H$%
\textsl{\ (cf. Lemma 1).} \textsl{The induced cohomology ring map }$\varphi
_{\beta _{1},\cdots ,\beta _{m};H}^{\ast }$\textsl{\ is given by}

\begin{center}
$\varphi _{\beta _{1},\cdots ,\beta _{m};H}^{\ast }(P_{w}(H))=(-1)^{r}%
\underset{_{\substack{ J=[i_{1},\cdots ,i_{r}]\subseteq \lbrack 1,\cdots ,m]
\\ \sigma _{J}=w}}}{\Sigma }x_{i_{1}}\cdots x_{i_{r}}$\textsl{,}
\end{center}

\noindent \textsl{where }$r=l(w)$.

\textbf{Proof.} Since $T\subseteq H\subset G$, the map $\varphi _{\beta
_{1},\cdots ,\beta _{m};H}$ factors through $\varphi _{\beta _{1},\cdots
,\beta _{m};T}$ in the fashion

\begin{center}
$%
\begin{array}{ccc}
\Gamma (\beta _{1},\cdots ,\beta _{m}) & \overset{\varphi _{\beta
_{1},\cdots ,\beta _{m};T}}{\rightarrow } & G/T \\
& \searrow & \downarrow \pi \\
\varphi _{\beta _{1},\cdots ,\beta _{m};H} &  & G/H%
\end{array}%
$,
\end{center}

\noindent where $\pi $ is the standard fibration with fiber $H/T$. By [D$%
_{3} $,Lemma 5.1] we have

\begin{quote}
(a)\textsl{\ Lemma 5 holds for the case }$H=T$\textsl{.}
\end{quote}

\noindent From [BGG, \S 5] we find that

\begin{quote}
(b)\textsl{\ the induced map }$\pi ^{\ast }:H^{\ast }(G/H)\rightarrow
H^{\ast }(G/T)$\textsl{\ is given by}
\end{quote}

\begin{center}
$\pi ^{\ast }(P_{w}(H))=P_{w}(T)$\textsl{,\ }$w\in \overline{W}\subset W$%
\textsl{.}
\end{center}

\noindent Combining (a) and (b) verifies Lemma 5.$\square $

\begin{center}
\textbf{4. Proof of the Theorem.}
\end{center}

For a $u\in \overline{W}$ with $l(u)=r$ and a $k\geq 1$, we assume as in
(1.2) that

\begin{enumerate}
\item[(3.1)] \noindent $\qquad \qquad \ \ \ \mathcal{P}^{k}(P_{u}(H))\equiv
\underset{l(v)=r+p(k-1),v\in \overline{W}}{\Sigma }a_{v,u}^{k}P_{v}(H)$, $%
a_{v,u}^{k}\in \mathbb{Z}_{p}$.
\end{enumerate}

\noindent Let $w=\sigma _{\beta _{1}}\circ \cdots \circ \sigma _{\beta _{m}}$%
, $\beta _{i}\in \Delta $, be a reduced decomposition of a $w\in \overline{W}
$ with $m=r+k(p-1)$, and let $A_{w}=(a_{i,j})_{m\times m}$ be the associated
Cartan matrix.

Let $\varphi _{\beta _{1},\cdots ,\beta _{m};H}:$ $\Gamma (\beta _{1},\cdots
,\beta _{m})\rightarrow G/H$ be the K-cycle associated to the ordered
sequence $(\beta _{1},\cdots ,\beta _{m})$ of simple roots. Applying the
ring map $\varphi _{\beta _{1},\cdots ,\beta _{m};H}^{\ast }$ to the
equation (3.1) in $H^{\ast }(G/H;\mathbb{Z}_{p})$ yields in $H^{\ast
}(\Gamma (\beta _{1},\cdots ,\beta _{m});\mathbb{Z}_{p})$ that

\begin{center}
$\varphi _{\beta _{1},\cdots ,\beta _{m};H}^{\ast }\mathcal{P}%
^{k}(P_{u}(H))\equiv \varphi _{\beta _{1},\cdots ,\beta _{m};H}^{\ast }[%
\underset{l(v)=r+p(k-1)}{\sum }a_{v,u}^{k}P_{v}(H)]$

$\equiv (-1)^{m}a_{w,u}^{k}x_{1}\cdots x_{m}$,
\end{center}

\noindent where the second equality follows from

\begin{center}
$\varphi _{\beta _{1},\cdots ,\beta _{m};H}^{\ast }[P_{v}(H)]=\{%
\begin{array}{c}
(-1)^{m}x_{1}\cdots x_{m}\text{ if }v=w\text{;} \\
0\text{ if }v\neq w\qquad \qquad%
\end{array}%
$.
\end{center}

\noindent by Lemma 5. On the other hand

\begin{center}
$\varphi _{\beta _{1},\cdots ,\beta _{m};H}^{\ast }\mathcal{P}%
^{k}(P_{u}(H))\equiv \mathcal{P}^{k}\{\varphi _{\beta _{1},\cdots ,\beta
_{m};H}^{\ast }(P_{u}(H))\}$

$\equiv (-1)^{r}\mathcal{P}^{k}\{\underset{_{\substack{ J=[i_{1},\cdots
,i_{r}]\subseteq \lbrack 1,\cdots ,m]  \\ \sigma _{J}=u}}}{\sum }%
x_{i_{1}}\cdots x_{i_{r}}\}$ \qquad (by Lemma 5)\quad \qquad $\equiv (-1)^{r}%
\underset{_{\substack{ J=[i_{1},\cdots ,i_{r}]\subseteq \lbrack 1,\cdots ,m]
\\ \sigma _{J}=u}}}{\sum }m_{k,p}(x_{i_{1}},\cdots ,x_{i_{r}})$ \quad (by
Lemma 4),
\end{center}

\noindent where the first equality comes from the naturality of $\mathcal{P}%
^{k}$ [SE]. Summarizing, we get in $H^{2m}(\Gamma (\beta _{1},\cdots ,\beta
_{m});\mathbb{Z}_{p})=\mathbb{Z}_{p}$ that

\begin{center}
$\qquad \underset{_{\substack{ J=[i_{1},\cdots ,i_{r}]\subseteq \lbrack
1,\cdots ,m]  \\ \sigma _{J}=u}}}{\sum }m_{k,p}(x_{i_{1}},\cdots
,x_{i_{r}})\equiv (-1)^{k(p-1)}a_{w,u}^{k}x_{1}\cdots x_{m}$.
\end{center}

\noindent Evaluating both sides on the orientation class $[\Gamma (\beta
_{1},\cdots ,\beta _{m})]$ \textsl{\ mod p} and noting that $p$ is a prime,
we get from Lemma 3 that

\begin{center}
$a_{w,u}^{k}\equiv <\underset{_{\substack{ J=[i_{1},\cdots ,i_{r}]\subseteq
\lbrack 1,\cdots ,m]  \\ \sigma _{J}=u}}}{\sum }m_{k,p}(x_{i_{1}},\cdots
,x_{i_{r}}),[\Gamma (\beta _{1},\cdots ,\beta _{m})]>$

$\equiv T_{A_{w}}(\underset{_{\substack{ J=[i_{1},\cdots ,i_{r}]\subseteq
\lbrack 1,\cdots ,m]  \\ \sigma _{J}=u}}}{\sum }m_{k,p}(x_{i_{1}},\cdots
,x_{i_{r}}))$.\qquad \qquad \qquad \qquad
\end{center}

\noindent This completes the proof.

\begin{center}
\textbf{5. Applications}
\end{center}

The Theorem handles the Problem C in its natural generality in the sense
that it applies uniformly to

\begin{quote}
(1) every flag manifold $G/H$;

(2) each Schubert classe in a given $G/H$; and

(3) is valid for every $k\geq 1$ and a prime $p\geq 2$.
\end{quote}

\noindent Because of these reasons a single program can be composed to
perform computation in various $G/H$ (cf. the discussion at the end of
Section 2). We list computational results from the program for some cases of
$G/H$.

For the Lie groups $G$ concerned below, we let a set of simple roots $\Delta
=\{\beta _{1},\cdots ,\beta _{n}\}$ of $G$ be given and ordered as that in
[Hu, p.64-75].

\bigskip

For a centralizer $H$ of a one-parameter subgroup in $G$, write $\overline{W}%
^{r}$ for the subset of $\overline{W}=W/W^{^{\prime }}$ consisting of the
elements with length $r$ (cf. Definition 1), where $W$ (resp. $W^{\prime }$)
is the Weyl group of $G$ (resp. $H$). The set $\{P_{w}(H)\mid $\textsl{\ }$%
w\in \overline{W}^{r}\}$\textsl{\ }forms a basis for the $2r$-dimensional
cohomology $H^{2r}(G/H)$ by Lemma 1.

Let $\Delta _{H}\subset \Delta $ be the subset consisting of simple roots of
$H$. Starting from the order on $\Delta $ as well as the subset $\Delta
_{H}\subset \Delta $, programs to decompose each $w\in \overline{W}^{r}$
uniquely into a reduced product, called \textsl{the minimal reduced
decomposition of }$w$, have been composed (cf. [DZZ,DZ]). If two $w$, $%
w^{\prime }\in \overline{W}^{r}$ are given by their minimal reduced
decompositions

\begin{center}
$w=\sigma _{\beta _{i_{1}}}\circ \cdots \circ \sigma _{\beta _{i_{r}}}$, $%
w^{\prime }=\sigma _{\beta _{j_{1}}}\circ \cdots \circ \sigma _{\beta
_{j_{r}}}$,
\end{center}

\noindent we say $w<w^{\prime }$ if there exists a $1\leq d<r$ such that $%
(i_{1},\cdots ,i_{d})=(j_{1},\cdots ,j_{d})$, but $i_{d+1}<j_{d+1}$. With
respect to this order $\overline{W}^{r}$ becomes an ordered set, hence can
be written as $\overline{W}^{r}=\{w_{r,i}\mid 1\leq i\leq \#\overline{W}%
^{r}\}$.

\bigskip

We write $\sigma _{i}$ instead of $\sigma _{\beta _{i}}$, $\beta _{i}\in
\Delta $. We list in Table A all $w_{r,i}\in \overline{W}$ by their minimal
reduced decompositions; followed by Table B that expresses all nontrivial $%
\mathcal{P}^{k}(s_{r,i})$, where the notion $s_{r,i}$ is used to simplify
the Schubert class $P_{w_{r,i}}(H)$ \textsl{\ mod p}.

\bigskip

\textbf{Example 1.} $G=G_{2}$ (the exceptional group of rank 2) and $H=T$ (a
maximal torus).

{\setlength{\tabcolsep}{1 mm} }

\begin{center}
\begin{tabular}{|l|l||l|l||l|l|}
\multicolumn{6}{l}{\textbf{Table A.} \textsl{\ Elements of }$\overline{W}$%
\textsl{\ and their minimal reduced decompositions.}} \\ \hline
${w}_{r,i}$ & decomposition & ${w}_{r,i}$ & decomposition & ${w}_{r,i}$ &
decomposition \\ \hline
${w}_{1,1}$ & $\sigma _{1}$ & ${w}_{3,1}$ & $\sigma _{1}\sigma _{2}\sigma
_{1}$ & ${w}_{5,1}$\ \  & $\sigma _{1}\sigma _{2}\sigma _{1}\sigma
_{2}\sigma _{1}$ \\ \hline
${w}_{1,2}$ & $\sigma _{2}$ & ${w}_{3,2}$ & $\sigma _{2}\sigma _{1}\sigma
_{2}$ & ${w}_{5,2}$\ \  & $\sigma _{2}\sigma _{1}\sigma _{2}\sigma
_{1}\sigma _{2}$ \\ \hline
${w}_{2,1}$ & $\sigma _{1}\sigma _{2}$ & ${w}_{4,1}$\ \  & $\sigma
_{1}\sigma _{2}\sigma _{1}\sigma _{2}$ & ${w}_{6,1}$\ \  & $\sigma
_{1}\sigma _{2}\sigma _{1}\sigma _{2}\sigma _{1}\sigma _{2}$ \\ \hline
${w}_{2,2}$ & $\sigma _{2}\sigma _{1}$ & ${w}_{4,2}$\ \  & $\sigma
_{2}\sigma _{1}\sigma _{2}\sigma _{1}$ &  &  \\ \hline
\end{tabular}

\bigskip

\begin{tabular}{|l|l|l|}
\multicolumn{3}{l}{\textbf{Table B.} \textsl{\ Nontrivial }$P^{k}(s_{r,i})$}
\\ \hline
${s}_{r,i}$ & $\mathcal{P}^{1}(s_{r,i})$ ($p=3$) & $\mathcal{P}^{1}(s_{r,i})$
($p=5$) \\ \hline
${s}_{1,1}$ & 0 & $3\,{s}_{5,1}$ \\ \hline
${s}_{1,2}$ & $2\,{s}_{3,2}$ & $2\,{s}_{5,2}$ \\ \hline
${s}_{2,1}$ & ${s}_{4,1}$ & 0 \\ \hline
${s}_{2,2}$ & $2\,{s}_{4,2}$ & 0 \\ \hline
${s}_{3,1}$ & ${s}_{5,1}$ & 0 \\ \hline
\end{tabular}

\bigskip
\end{center}

\textbf{Example 2.} $G=F_{4}$ (the exceptional group of rank 4) and $%
H=Spin(7)\times S^{1}$.

{\setlength{\tabcolsep}{0mm} }

\begin{center}
\begin{tabular}{cc}
\multicolumn{2}{c}{\textbf{Table A.} \textsl{\ Elements of }$\overline{W}$}
\\
\multicolumn{2}{c}{\textsl{\ and their minimal reduced decompositions.}} \\
\begin{tabular}{|l|l|}
\hline
\ ${w}_{r,i}$ \  & \ decomposition \\ \hline
\ ${w}_{1,1}$\  & \ $\sigma_1$ \\ \hline
\ ${w}_{2,1}$\  & \ $\sigma_2\sigma_1$ \\ \hline
\ ${w}_{3,1}$\  & \ $\sigma_3\sigma_2\sigma_1$ \\ \hline
\ ${w}_{4,1}$\  & \ $\sigma_2\sigma_3\sigma_2\sigma_1$ \\ \hline
\ ${w}_{4,2}$\  & \ $\sigma_4\sigma_3\sigma_2\sigma_1$ \\ \hline
\ ${w}_{5,1}$\  & \ $\sigma_1\sigma_2\sigma_3\sigma_2\sigma_1$ \\ \hline
\ ${w}_{5,2}$\  & \ $\sigma_2\sigma_4\sigma_3\sigma_2\sigma_1$ \\ \hline
\ ${w}_{6,1}$\  & \ $\sigma_1\sigma_2\sigma_4\sigma_3\sigma_2\sigma_1 $ \\
\hline
\ ${w}_{6,2}$\  & \ $\sigma_3\sigma_2\sigma_4\sigma_3\sigma_2\sigma_1 $ \\
\hline
\ ${w}_{7,1}$\  & \ $\sigma_1\sigma_3\sigma_2\sigma_4\sigma_3\sigma_2
\sigma_1$ \\ \hline
\ ${w}_{7,2}$\  & \ $\sigma_2\sigma_3\sigma_2\sigma_4\sigma_3\sigma_2
\sigma_1$ \\ \hline
\ ${w}_{8,1}$\ \  & \ $\sigma_1\sigma_2\sigma_3\sigma_2\sigma_4\sigma_3
\sigma_2\sigma_1$\ \  \\ \hline
\end{tabular}
&
\begin{tabular}{|l|l|}
\hline
\ ${w}_{r,i}$ & \ decomposition \\ \hline
\ ${w}_{8,2}$\  & \ $\sigma_2\sigma_1\sigma_3\sigma_2\sigma_4\sigma_3
\sigma_2\sigma_1$ \\ \hline
\ ${w}_{9,1}$\  & \ $\sigma_1\sigma_2\sigma_1\sigma_3\sigma_2\sigma_4
\sigma_3\sigma_2\sigma_1$ \\ \hline
\ ${w}_{9,2}$\  & \ $\sigma_3\sigma_2\sigma_1\sigma_3\sigma_2\sigma_4
\sigma_3\sigma_2\sigma_1$ \\ \hline
\ ${w}_{10,1}$\  & \ $\sigma_1\sigma_3\sigma_2\sigma_1\sigma_3\sigma_2
\sigma_4\sigma_3\sigma_2\sigma_1$ \\ \hline
\ ${w}_{10,2}$\  & \ $\sigma_4\sigma_3\sigma_2\sigma_1\sigma_3\sigma_2
\sigma_4\sigma_3\sigma_2\sigma_1$ \\ \hline
\ ${w}_{11,1}$\  & \ $\sigma_1\sigma_4\sigma_3\sigma_2\sigma_1\sigma_3
\sigma_2\sigma_4\sigma_3\sigma_2\sigma_1$ \\ \hline
\ ${w}_{11,2}$\  & \ $\sigma_2\sigma_1\sigma_3\sigma_2\sigma_1\sigma_3
\sigma_2\sigma_4\sigma_3\sigma_2\sigma_1$ \\ \hline
\ ${w}_{12,1}$\  & \ $\sigma_2\sigma_1\sigma_4\sigma_3\sigma_2\sigma_1
\sigma_3\sigma_2\sigma_4\sigma_3\sigma_2\sigma_1 $ \\ \hline
\ ${w}_{13,1}$\  & \ $\sigma_3\sigma_2\sigma_1\sigma_4\sigma_3\sigma_2
\sigma_1\sigma_3\sigma_2\sigma_4\sigma_3\sigma_2 \sigma_1$ \\ \hline
\ ${w}_{14,1}$\  & \ $\sigma_2\sigma_3\sigma_2\sigma_1\sigma_4\sigma_3
\sigma_2\sigma_1\sigma_3\sigma_2\sigma_4\sigma_3 \sigma_2\sigma_1$ \\ \hline
\ ${w}_{15,1}$\ \  & \ $\sigma_1\sigma_2\sigma_3\sigma_2\sigma_1\sigma_4
\sigma_3\sigma_2\sigma_1\sigma_3\sigma_2\sigma_4 \sigma_3\sigma_2\sigma_1$\
\  \\ \hline
&  \\ \hline
\end{tabular}%
\end{tabular}
\end{center}

{\setlength{\tabcolsep}{1mm} \lineskip0 mm }

\begin{center}
{\footnotesize \setlength{\tabcolsep}{0.8mm}
\begin{tabular}{|l|l|l|l|l|l|l|l|}
\multicolumn{8}{c}{\textbf{Table B.} \textsl{\ Nontrivial }$P^{k}(s_{r,i})$}
\\ \hline
${s}_{r,i}$ & \multicolumn{2}{c|}{$p=3$} & \multicolumn{3}{c|}{$p=5$} &
\multicolumn{2}{c|}{$p=7$} \\ \cline{2-8}
& $\mathcal{P}^{1}(s_{r,i})$ & $\mathcal{P}^{2}(s_{r,i})$ & $\mathcal{P}%
^{1}(s_{r,i})$ & $\mathcal{P}^{2}(s_{r,i})$ & $\mathcal{P}^{3}(s_{r,i})$ & $%
\mathcal{P}^{1}(s_{r,i})$ & $\mathcal{P}^{2}(s_{r,i})$ \\ \hline
${s}_{1,1}$ & ${s}_{3,1}$ & 0 & ${s}_{5,1} + 2\, {s}_{5,2}$ & 0 & 0 & $5\, {s%
}_{7,1} + 2\, {s}_{7,2}$ & 0 \\ \hline
${s}_{2,1}$ & $2\, {s}_{4,1} + 2\, {s}_{4,2}$ & $2\, {s}_{6,2}$ & ${s}_{6,1}
+ 4\, {s}_{6,2}$ & {$%
\begin{array}{ll}
3 {s}_{10,1} &  \\
+ 2 {s}_{10,2} &
\end{array}%
$} & 0 & $4\, {s}_{8,1} + 3\, {s}_{8,2}$ & ${s}_{14,1}$ \\ \hline
${s}_{3,1}$ & 0 & $0$ & ${s}_{7,2}$ & $4\, {s}_{11,2}$ & $3\, {s}_{15,1}$ & $%
{s}_{9,2}$ & $3\, {s}_{15,1}$ \\ \hline
${s}_{4,1}$ & ${s}_{6,2}$ & 0 & $2\, {s}_{8,1} + 4\, {s}_{8,2}$ & 0 &  & $%
3\, {s}_{10,1} + 3\, {s}_{10,2}$ &  \\ \hline
${s}_{4,2}$ & ${s}_{6,2}$ & 0 & $4\, {s}_{8,1} + 4\, {s}_{8,2}$ & 0 &  & $%
3\, {s}_{10,1} + 3\, {s}_{10,2}$ &  \\ \hline
${s}_{5,1}$ & ${s}_{7,1}$ & 0 & ${s}_{9,1} + 2\, {s}_{9,2}$ & 0 &  & $5\, {s}%
_{11,2}$ &  \\ \hline
${s}_{5,2}$ & $2\, {s}_{7,2}$ & 0 & $2\, {s}_{9,1} + 4\, {s}_{9,2}$ & 0 &  &
$4\, {s}_{11,1} + 3\, {s}_{11,2}$ &  \\ \hline
${s}_{6,1}$ & $2\, {s}_{8,1} + {s}_{8,2}$ & $%
\begin{array}{ll}
2 {s}_{10,1} &  \\
+ 2 {s}_{10,2} &
\end{array}%
$ & $3 {s}_{10,1} + 4 {s}_{10,2}$ & 0 &  & $3\, {s}_{12,1}$ &  \\ \hline
${s}_{6,2}$ & 0 & 0 & $2\, {s}_{10,1}$ & 0 &  & $2\, {s}_{12,1}$ &  \\ \hline
${s}_{7,1}$ & 0 & 0 & $4 {s}_{11,1} + 3 {s}_{11,2}$ & 0 &  & ${s}_{13,1}$ &
\\ \hline
${s}_{7,2}$ & 0 & 0 & $3\, {s}_{11,2}$ & 0 &  & ${s}_{13,1}$ &  \\ \hline
${s}_{8,2}$ & ${s}_{10,1} + {s}_{10,2}$ & 0 & $2\, {s}_{12,1}$ &  &  & $3\, {%
s}_{14,1}$ &  \\ \hline
${s}_{9,1}$ & ${s}_{11,1} + {s}_{11,2}$ & 0 & $2\, {s}_{13,1}$ &  &  & $2\, {%
s}_{15,1}$ &  \\ \hline
${s}_{9,2}$ & ${s}_{11,1} + 2\, {s}_{11,2}$ & 0 & $0$ &  &  & $6\, {s}%
_{15,1} $ &  \\ \hline
${s}_{10,1}$ & ${s}_{12,1}$ & $2\, {s}_{14,1}$ & $4\, {s}_{14,1}$ &  &  &  &
\\ \hline
${s}_{10,2}$ & $2\, {s}_{12,1}$ & ${s}_{14,1}$ & $4\, {s}_{14,1}$ &  &  &  &
\\ \hline
${s}_{11,2}$ & 0 & 0 & ${s}_{15,1}$ &  &  &  &  \\ \hline
${s}_{12,1}$ & ${s}_{14,1}$ &  &  &  &  &  &  \\ \hline
${s}_{13,1}$ & $2\, {s}_{15,1}$ &  &  &  &  &  &  \\ \hline
\end{tabular}
}
\end{center}

\textbf{Example 3.} $G=SO(12)$ (the special orthogonal group of order $12$)
and $H=U(6)$. The flag manifold $G/H$ is the Grassmannian of complex
structures on the $12$-dimensional real Euclidean space $\mathbb{R}^{12}$ [D$%
_{2}$].

{\small \setlength{\tabcolsep}{0mm} }

\begin{center}
{\small \setlength{\tabcolsep}{0 mm}
\begin{tabular}{cc}
\multicolumn{2}{c}{\textbf{Table A.} \textsl{\ Elements of }$\overline{W}$
and their minimal reduced decompositions.} \\
\begin{tabular}{|l|l|}
\hline
\ ${w}_{r,i}$ & \ decomposition \\ \hline
\ ${w}_{1,1}$ & \ $\sigma_6$ \\ \hline
\ ${w}_{2,1}$ & \ $\sigma_4\sigma_6$ \\ \hline
\ ${w}_{3,1}$ & \ $\sigma_3\sigma_4\sigma_6$ \\ \hline
\ ${w}_{3,2}$ & \ $\sigma_5\sigma_4\sigma_6$ \\ \hline
\ ${w}_{4,1}$ & \ $\sigma_2\sigma_3\sigma_4\sigma_6$ \\ \hline
\ ${w}_{4,2}$ & \ $\sigma_3\sigma_5\sigma_4\sigma_6$ \\ \hline
\ ${w}_{5,1}$ & \ $\sigma_1\sigma_2\sigma_3\sigma_4\sigma_6$ \\ \hline
\ ${w}_{5,2}$ & \ $\sigma_2\sigma_3\sigma_5\sigma_4\sigma_6$ \\ \hline
\ ${w}_{5,3}$ & \ $\sigma_4\sigma_3\sigma_5\sigma_4\sigma_6$ \\ \hline
\ ${w}_{6,1}$ & \ $\sigma_1\sigma_2\sigma_3\sigma_5\sigma_4\sigma_6 $ \\
\hline
\ ${w}_{6,2}$ & \ $\sigma_2\sigma_4\sigma_3\sigma_5\sigma_4\sigma_6 $ \\
\hline
\ ${w}_{6,3}$ & \ $\sigma_6\sigma_4\sigma_3\sigma_5\sigma_4\sigma_6 $ \\
\hline
\ ${w}_{7,1}$ & \ $\sigma_1\sigma_2\sigma_4\sigma_3\sigma_5\sigma_4 \sigma_6$
\\ \hline
\ ${w}_{7,2}$ & \ $\sigma_2\sigma_6\sigma_4\sigma_3\sigma_5\sigma_4 \sigma_6$
\\ \hline
\ ${w}_{7,3}$ & \ $\sigma_3\sigma_2\sigma_4\sigma_3\sigma_5\sigma_4 \sigma_6$
\\ \hline
\ ${w}_{8,1}\ \ $ & \ $\sigma_1\sigma_2\sigma_6\sigma_4\sigma_3\sigma_5
\sigma_4\sigma_6$ \ \  \\ \hline
\end{tabular}
&
\begin{tabular}{|l|l|}
\hline
\ ${w}_{r,i}$ & \ decomposition \\ \hline
\ ${w}_{8,2}$ & \ $\sigma_1\sigma_3\sigma_2\sigma_4\sigma_3\sigma_5
\sigma_4\sigma_6$ \\ \hline
\ ${w}_{8,3}$ & \ $\sigma_3\sigma_2\sigma_6\sigma_4\sigma_3\sigma_5
\sigma_4\sigma_6$ \\ \hline
\ ${w}_{9,1}$ & \ $\sigma_1\sigma_3\sigma_2\sigma_6\sigma_4\sigma_3
\sigma_5\sigma_4\sigma_6$ \\ \hline
\ ${w}_{9,2}$ & \ $\sigma_2\sigma_1\sigma_3\sigma_2\sigma_4\sigma_3
\sigma_5\sigma_4\sigma_6$ \\ \hline
\ ${w}_{9,3}$ & \ $\sigma_4\sigma_3\sigma_2\sigma_6\sigma_4\sigma_3
\sigma_5\sigma_4\sigma_6$ \\ \hline
\ ${w}_{10,1}$ & \ $\sigma_1\sigma_4\sigma_3\sigma_2\sigma_6\sigma_4
\sigma_3\sigma_5\sigma_4\sigma_6$ \\ \hline
\ ${w}_{10,2}$ & \ $\sigma_2\sigma_1\sigma_3\sigma_2\sigma_6\sigma_4
\sigma_3\sigma_5\sigma_4\sigma_6$ \\ \hline
\ ${w}_{10,3}$ & \ $\sigma_5\sigma_4\sigma_3\sigma_2\sigma_6\sigma_4
\sigma_3\sigma_5\sigma_4\sigma_6$ \\ \hline
\ ${w}_{11,1}$ & \ $\sigma_1\sigma_5\sigma_4\sigma_3\sigma_2\sigma_6
\sigma_4\sigma_3\sigma_5\sigma_4\sigma_6$ \\ \hline
\ ${w}_{11,2}$ & \ $\sigma_2\sigma_1\sigma_4\sigma_3\sigma_2\sigma_6
\sigma_4\sigma_3\sigma_5\sigma_4\sigma_6$ \\ \hline
\ ${w}_{12,1}$ & \ $\sigma_2\sigma_1\sigma_5\sigma_4\sigma_3\sigma_2
\sigma_6\sigma_4\sigma_3\sigma_5\sigma_4\sigma_6 $ \\ \hline
\ ${w}_{12,2}$ & \ $\sigma_3\sigma_2\sigma_1\sigma_4\sigma_3\sigma_2
\sigma_6\sigma_4\sigma_3\sigma_5\sigma_4\sigma_6 $ \\ \hline
\ ${w}_{13,1}$ & \ $\sigma_3\sigma_2\sigma_1\sigma_5\sigma_4\sigma_3
\sigma_2\sigma_6\sigma_4\sigma_3\sigma_5\sigma_4 \sigma_6$ \\ \hline
\ ${w}_{14,1}$ & \ $\sigma_4\sigma_3\sigma_2\sigma_1\sigma_5\sigma_4
\sigma_3\sigma_2\sigma_6\sigma_4\sigma_3\sigma_5 \sigma_4\sigma_6$ \\ \hline
\ ${w}_{15,1}$\ \  & \ $\sigma_6\sigma_4\sigma_3\sigma_2\sigma_1\sigma_5
\sigma_4\sigma_3\sigma_2\sigma_6\sigma_4\sigma_3 \sigma_5\sigma_4\sigma_6$\
\  \\ \hline
&  \\ \hline
\end{tabular}%
\end{tabular}
}
\end{center}

{\setlength{\tabcolsep}{1mm} }

\begin{center}
{\footnotesize
\begin{tabular}{|l|l|l|l|l|l|}
\multicolumn{6}{c}{\textbf{Table B${}_1$.} \textsl{\ Nontrivial }$%
P^{k}(s_{r,i})$ \textsl{\ for $p=3$}} \\ \hline
${s}_{r,i}$ & $\mathcal{P}^{1}(s_{r,i})$ & $\mathcal{P}^{2}(s_{r,i})$ & $%
\mathcal{P}^{3}(s_{r,i})$ & $\mathcal{P}^{4}(s_{r,i})$ & $\mathcal{P}%
^{5}(s_{r,i})$ \\ \hline
${s}_{1,1}$ & ${s}_{3,1} + {s}_{3,2}$ & $0$ & $0$ & $0$ & $0$ \\ \hline
${s}_{2,1}$ & $2\, {s}_{4,1} + {s}_{4,2}$ & ${s}_{6,1} + 2\, {s}_{6,2} + 2\,
{s}_{6,3}$ & $0$ & $0$ & $0$ \\ \hline
${s}_{3,1}$ & ${s}_{5,2}$ & $2\, {s}_{7,1} + 2\, {s}_{7,2} +2\, {s}_{7,3}$ &
${s}_{9,1} + {s}_{9,2} + 2 {s}_{9,3}$ & $0$ & $0$ \\ \hline
${s}_{3,2}$ & $2\, {s}_{5,2}$ & ${s}_{7,1} + {s}_{7,2} + {s}_{7,3}$ & $2 {s}%
_{9,1} + {s}_{9,2} + {s}_{9,3}$ & $0$ & $0$ \\ \hline
${s}_{4,1}$ & ${s}_{6,1}$ & $2\, {s}_{8,1} + 2\, {s}_{8,2}$ & $2\, {s}%
_{10,1} + {s}_{10,2}$ & $2\, {s}_{12,2}$ & $0$ \\ \hline
${s}_{4,2}$ & ${s}_{6,2} + {s}_{6,3}$ & $2\, {s}_{8,1} + 2\, {s}_{8,2}$ & $%
2\, {s}_{10,1} + 2\, {s}_{10,2}$ & ${s}_{12,1}$ & $0$ \\ \hline
${s}_{5,2}$ & ${s}_{7,1} + {s}_{7,2} + {s}_{7,3}$ & $0$ & ${s}_{11,2}$ & ${s}%
_{13,1}$ & $2\, {s}_{15,1}$ \\ \hline
${s}_{5,3}$ & ${s}_{7,2} + 2\, {s}_{7,3}$ & ${s}_{9,2}$ & ${s}_{11,2}$ & ${s}%
_{13,1}$ & $2\, {s}_{15,1}$ \\ \hline
${s}_{6,1}$ & ${s}_{8,1} + {s}_{8,2}$ & $0$ & $0$ & $0$ &  \\ \hline
${s}_{6,2}$ & ${s}_{8,1} + {s}_{8,2} + {s}_{8,3}$ & $2 {s}_{10,1} + 2 {s}%
_{10,2} + 2 {s}_{10,3}$ & ${s}_{12,1} + {s}_{12,2}$ & $0$ &  \\ \hline
${s}_{6,3}$ & $2\, {s}_{8,3}$ & ${s}_{10,1} + {s}_{10,2} + {s}_{10,3}$ & ${s}%
_{12,1} + {s}_{12,2}$ & $0$ &  \\ \hline
${s}_{7,1}$ & ${s}_{9,1} + 2\, {s}_{9,2}$ & $2\, {s}_{11,1} + 2\, {s}_{11,2}$
& $0$ & $0$ &  \\ \hline
${s}_{7,2}$ & ${s}_{9,1}$ & $2\, {s}_{11,1} + 2\, {s}_{11,2}$ & $2\, {s}%
_{13,1}$ & $2\, {s}_{15,1}$ &  \\ \hline
${s}_{7,3}$ & ${s}_{9,1} + {s}_{9,2}$ & $2\, {s}_{11,1} + 2\, {s}_{11,2}$ & $%
0$ &  &  \\ \hline
${s}_{8,1}$ & $2\, {s}_{10,2}$ & ${s}_{12,1} + {s}_{12,2}$ & ${s}_{14,1}$ &
&  \\ \hline
${s}_{8,2}$ & ${s}_{10,2}$ & $2\, {s}_{12,1} + 2\, {s}_{12,2}$ & $2\, {s}%
_{14,1}$ &  &  \\ \hline
${s}_{8,3}$ & ${s}_{10,1} + {s}_{10,2} + {s}_{10,3}$ & $0$ & ${s}_{14,1}$ &
&  \\ \hline
${s}_{9,1}$ & ${s}_{11,1} + {s}_{11,2}$ & $0$ & $0$ &  &  \\ \hline
${s}_{9,3}$ & ${s}_{11,1} + {s}_{11,2}$ & $0$ & ${s}_{15,1}$ &  &  \\ \hline
${s}_{10,1}$ & ${s}_{12,1} + 2\, {s}_{12,2}$ & $0$ &  &  &  \\ \hline
${s}_{10,2}$ & ${s}_{12,1} + {s}_{12,2}$ & $0$ &  &  &  \\ \hline
${s}_{10,3}$ & ${s}_{12,1}$ & $0$ &  &  &  \\ \hline
${s}_{11,1}$ & $2\, {s}_{13,1}$ & ${s}_{15,1}$ &  &  &  \\ \hline
${s}_{11,2}$ & ${s}_{13,1}$ & $2\, {s}_{15,1}$ &  &  &  \\ \hline
${s}_{13,1}$ & ${s}_{15,1}$ &  &  &  &  \\ \hline
\end{tabular}
}
\end{center}


\begin{center}
\begin{tabular}{|l|l|l|l|}
\multicolumn{4}{c}{\textbf{Table B${}_2$.} \textsl{\ Nontrivial }$%
P^{k}(s_{r,i})$ \textsl{for $p=5$}} \\ \hline
${s}_{r,i}$ & $\mathcal{P}^{1}(s_{r,i})$ & $\mathcal{P}^{2}(s_{r,i})$ & $%
\mathcal{P}^{3}(s_{r,i})$ \\ \hline
${s}_{1,1}$ & ${s}_{5,1} + 3\, {s}_{5,2} + 2\, {s}_{5,3}$ & $0$ & $0$ \\
\hline
${s}_{2,1}$ & $3\, {s}_{6,1} + 4\, {s}_{6,3}$ & $0$ & $0$ \\ \hline
${s}_{3,1}$ & $4\, {s}_{7,2} + 3\, {s}_{7,3}$ & $2\, {s}_{11,1}$ & $4\, {s}%
_{15,1}$ \\ \hline
${s}_{3,2}$ & $2\, {s}_{7,1} + 2\, {s}_{7,2} + 2\, {s}_{7,3}$ & ${s}_{11,1}$
& $2\, {s}_{15,1}$ \\ \hline
${s}_{4,1}$ & $4\, {s}_{8,1} + 3\, {s}_{8,2}$ & ${s}_{12,2}$ &  \\ \hline
${s}_{4,2}$ & $4\, {s}_{8,2} + 4\, {s}_{8,3}$ & $3\, {s}_{12,1} + 2\, {s}%
_{12,2}$ &  \\ \hline
${s}_{5,2}$ & $2\, {s}_{9,1} + 2\, {s}_{9,2} + \ 2\, {s}_{9,3}$ & $2\, {s}%
_{13,1}$ &  \\ \hline
${s}_{5,3}$ & $2\, {s}_{9,1} + 2\, {s}_{9,2} + \ 2\, {s}_{9,3}$ & $2\, {s}%
_{13,1}$ &  \\ \hline
${s}_{6,1}$ & $2\, {s}_{10,1} + 3\, {s}_{10,2}$ & $0$ &  \\ \hline
${s}_{6,2}$ & $2\, {s}_{10,1} + 2\, {s}_{10,2} + \ 4\, {s}_{10,3}$ & $4\, {s}%
_{14,1}$ &  \\ \hline
${s}_{6,3}$ & $3\, {s}_{10,1} + 2\, {s}_{10,2} + \ {s}_{10,3}$ & $0$ &  \\
\hline
${s}_{7,1}$ & $4\, {s}_{11,1}$ & $2\, {s}_{15,1}$ &  \\ \hline
${s}_{7,2}$ & $3\, {s}_{11,1}$ & $4\, {s}_{15,1}$ &  \\ \hline
${s}_{7,3}$ & $4\, {s}_{11,1}$ & $2\, {s}_{15,1}$ &  \\ \hline
${s}_{8,1}$ & $2\, {s}_{12,1} + 2\, {s}_{12,2}$ &  &  \\ \hline
${s}_{8,2}$ & $4\, {s}_{12,1} + 3\, {s}_{12,2}$ &  &  \\ \hline
${s}_{8,3}$ & $3\, {s}_{12,2}$ &  &  \\ \hline
${s}_{9,1}$ & $4\, {s}_{13,1}$ &  &  \\ \hline
${s}_{9,3}$ & $3\, {s}_{13,1}$ &  &  \\ \hline
${s}_{10,1}$ & $2\, {s}_{14,1}$ &  &  \\ \hline
${s}_{10,2}$ & $2\, {s}_{14,1}$ &  &  \\ \hline
${s}_{11,1}$ & ${s}_{15,1}$ &  &  \\ \hline
${s}_{11,2}$ & $4\, {s}_{15,1}$ &  &  \\ \hline
\end{tabular}
\end{center}


\textbf{Example 4.} $G=U(7)$ (the unitary group of order $7$) and $%
H=U(3)\times U(4)$. The flag manifold $G/H$ is the Grassmannian of $3$%
-planes through the origin in $\mathbb{C}^{7}$.

{\small \setlength{\tabcolsep}{0mm} \lineskip 0mm }

\begin{center}
{\small \setlength{\tabcolsep}{0mm} \lineskip 0mm
\begin{tabular}{c||c||c}
\multicolumn{3}{c}{\textbf{Table A.} \textsl{\ Elements of }$\overline{W}$%
\textsl{\ and their minimal reduced decompositions.}} \\
\begin{tabular}{|l|l}
\hline
\ $w_{r,i}$ & \ decomposition \\ \hline
\ ${w}_{1,1}$ & \ $\sigma_3$ \\ \hline
\ ${w}_{2,1}$ & \ $\sigma_2\sigma_3$ \\ \hline
\ ${w}_{2,2}$ & \ $\sigma_4\sigma_3$ \\ \hline
\ ${w}_{3,1}$ & \ $\sigma_1\sigma_2\sigma_3$ \\ \hline
\ ${w}_{3,2}$ & \ $\sigma_2\sigma_4\sigma_3$ \\ \hline
\ ${w}_{3,3}$ & \ $\sigma_5\sigma_4\sigma_3$ \\ \hline
\ ${w}_{4,1}$ & \ $\sigma_1\sigma_2\sigma_4\sigma_3$ \\ \hline
\ ${w}_{4,2}$ & \ $\sigma_2\sigma_5\sigma_4\sigma_3$ \\ \hline
\ ${w}_{4,3}$ & \ $\sigma_3\sigma_2\sigma_4\sigma_3$ \\ \hline
\ ${w}_{4,4}$ & \ $\sigma_6\sigma_5\sigma_4\sigma_3$ \\ \hline
\ ${w}_{5,1}$ & \ $\sigma_1\sigma_2\sigma_5\sigma_4\sigma_3$ \\ \hline
\ ${w}_{5,2}$ \  & \ $\sigma_1\sigma_3\sigma_2\sigma_4\sigma_3$\  \\ \hline
\end{tabular}
&
\begin{tabular}{|l|l|}
\hline
\ $w_{r,i}$ & \ decomposition \\ \hline
\ ${w}_{5,3}$ & \ $\sigma_2\sigma_6\sigma_5\sigma_4\sigma_3$ \\ \hline
\ ${w}_{5,4}$ & \ $\sigma_3\sigma_2\sigma_5\sigma_4\sigma_3$ \\ \hline
\ ${w}_{6,1}$ & \ $\sigma_1\sigma_2\sigma_6\sigma_5\sigma_4\sigma_3 $ \\
\hline
\ ${w}_{6,2}$ & \ $\sigma_1\sigma_3\sigma_2\sigma_5\sigma_4\sigma_3 $ \\
\hline
\ ${w}_{6,3}$ & \ $\sigma_2\sigma_1\sigma_3\sigma_2\sigma_4\sigma_3 $ \\
\hline
\ ${w}_{6,4}$ & \ $\sigma_3\sigma_2\sigma_6\sigma_5\sigma_4\sigma_3 $ \\
\hline
\ ${w}_{6,5}$ & \ $\sigma_4\sigma_3\sigma_2\sigma_5\sigma_4\sigma_3 $ \\
\hline
\ ${w}_{7,1}$ & \ $\sigma_1\sigma_3\sigma_2\sigma_6\sigma_5\sigma_4 \sigma_3$
\\ \hline
\ ${w}_{7,2}$ & \ $\sigma_1\sigma_4\sigma_3\sigma_2\sigma_5\sigma_4 \sigma_3$
\\ \hline
\ ${w}_{7,3}$ & \ $\sigma_2\sigma_1\sigma_3\sigma_2\sigma_5\sigma_4 \sigma_3$
\\ \hline
\ ${w}_{7,4}$ \  & \ $\sigma_4\sigma_3\sigma_2\sigma_6\sigma_5\sigma_4
\sigma_3$\ \  \\ \hline
&  \\ \hline
\end{tabular}
&
\begin{tabular}{l|l|}
\hline
\ $w_{r,i}$ & \ decomposition \\ \hline
\ ${w}_{8,1}$ & \ $\sigma_1\sigma_4\sigma_3\sigma_2\sigma_6\sigma_5
\sigma_4\sigma_3$ \\ \hline
\ ${w}_{8,2}$ & \ $\sigma_2\sigma_1\sigma_3\sigma_2\sigma_6\sigma_5
\sigma_4\sigma_3$ \\ \hline
\ ${w}_{8,3}$ & \ $\sigma_2\sigma_1\sigma_4\sigma_3\sigma_2\sigma_5
\sigma_4\sigma_3$ \\ \hline
\ ${w}_{8,4}$ & \ $\sigma_5\sigma_4\sigma_3\sigma_2\sigma_6\sigma_5
\sigma_4\sigma_3$ \\ \hline
\ ${w}_{9,1}$ & \ $\sigma_1\sigma_5\sigma_4\sigma_3\sigma_2\sigma_6
\sigma_5\sigma_4\sigma_3$ \\ \hline
\ ${w}_{9,2}$ & \ $\sigma_2\sigma_1\sigma_4\sigma_3\sigma_2\sigma_6
\sigma_5\sigma_4\sigma_3$ \\ \hline
\ ${w}_{9,3}$ & \ $\sigma_3\sigma_2\sigma_1\sigma_4\sigma_3\sigma_2
\sigma_5\sigma_4\sigma_3$ \\ \hline
\ ${w}_{10,1}$ & \ $\sigma_2\sigma_1\sigma_5\sigma_4\sigma_3\sigma_2
\sigma_6\sigma_5\sigma_4\sigma_3$ \\ \hline
\ ${w}_{10,2}$ & \ $\sigma_3\sigma_2\sigma_1\sigma_4\sigma_3\sigma_2
\sigma_6\sigma_5\sigma_4\sigma_3$ \\ \hline
\ ${w}_{11,1}$ & \ $\sigma_3\sigma_2\sigma_1\sigma_5\sigma_4\sigma_3
\sigma_2\sigma_6\sigma_5\sigma_4\sigma_3$ \\ \hline
\ ${w}_{12,1}$ \  & \ $\sigma_4\sigma_3\sigma_2\sigma_1\sigma_5\sigma_4
\sigma_3\sigma_2\sigma_6\sigma_5\sigma_4\sigma_3$\ \  \\ \hline
&  \\ \hline
\end{tabular}%
\end{tabular}
}
\end{center}

{\setlength{\tabcolsep}{1mm} }

\begin{center}
{\footnotesize \setlength{\tabcolsep}{1mm}
\begin{tabular}{|l|l|l|l|l|}
\multicolumn{5}{c}{\textbf{Table B${}_1$.} \textsl{\ Nontrivial }$%
P^{k}(s_{r,i})$ \textsl{for $p=3$}} \\ \hline
$s_{r,i}$ & $\mathcal{P}^{1}(s_{r,i})$ & $\mathcal{P}^{2}(s_{r,i})$ & $%
\mathcal{P}^{3}(s_{r,i})$ & $\mathcal{P}^{4}(s_{r,i})$ \\ \hline
${s}_{1,1}$ & ${s}_{3,1} + 2\, {s}_{3,2} + {s}_{3,3}$ & $0$ & $0$ & $0$ \\
\hline
${s}_{2,1}$ & $2\, {s}_{4,1} + {s}_{4,2} + 2\, {s}_{4,3}$ & $2\, {s}_{6,2} +
{s}_{6,3} + {s}_{6,5}$ & $0$ & $0$ \\ \hline
${s}_{2,2}$ & ${s}_{4,1} + 2 {s}_{4,2} + 2 {s}_{4,3} + 2 {s}_{4,4}$ & ${s}%
_{6,1} + 2\, {s}_{6,2} + {s}_{6,3} + \ {s}_{6,5}$ & $0$ & $0$ \\ \hline
${s}_{3,1}$ & ${s}_{5,1} + 2\, {s}_{5,2}$ & ${s}_{7,2} + 2\, {s}_{7,3}$ & ${s%
}_{9,3}$ & $0$ \\ \hline
${s}_{3,2}$ & $2 {s}_{5,1} + 2 {s}_{5,2} + 2 {s}_{5,3} + 2 {s}_{5,4}$ & $2\,
{s}_{7,1} + 2\, {s}_{7,4}$ & ${s}_{9,1} + 2\, {s}_{9,2} + 2\, {s}_{9,3}$ & $%
0 $ \\ \hline
${s}_{3,3}$ & ${s}_{5,1} + 2\, {s}_{5,3} + 2\, {s}_{5,4}$ & $2 {s}_{7,1} + 2
{s}_{7,2} + {s}_{7,3} + 2 {s}_{7,4}$ & ${s}_{9,1} + 2\, {s}_{9,2} + {s}%
_{9,3} $ & $0$ \\ \hline
${s}_{4,1}$ & $2\, {s}_{6,1} + 2\, {s}_{6,2}$ & $\, {s}_{8,1} + 2\, {s}%
_{8,2} + 2\, {s}_{8,3}$ & ${s}_{10,1} + {s}_{10,2}$ & ${s}_{12,1}$ \\ \hline
${s}_{4,2}$ & $2\, {s}_{6,1} + 2\, {s}_{6,2} + 2\, {s}_{6,4}$ & ${s}_{8,1} +
2\, {s}_{8,3} + {s}_{8,4}$ & $2\, {s}_{10,2}$ & $2 {s}_{12,1}$ \\ \hline
${s}_{4,3}$ & $2\, {s}_{6,2} + {s}_{6,3} + 2\, {s}_{6,4} \ + {s}_{6,5}$ & $%
2\, {s}_{8,1} + {s}_{8,2} + {s}_{8,4}$ & $2\, {s}_{10,1} + {s}_{10,2}$ & ${s}%
_{12,1}$ \\ \hline
${s}_{4,4}$ & ${s}_{6,1} + 2\, {s}_{6,4}$ & $2\, {s}_{8,1} + {s}_{8,2} + {s}%
_{8,4}$ & $2\, {s}_{10,1} + {s}_{10,2}$ & ${s}_{12,1}$ \\ \hline
${s}_{5,1}$ & $2\, {s}_{7,1}$ & ${s}_{9,1} + 2\, {s}_{9,2}$ & $0$ &  \\
\hline
${s}_{5,2}$ & $2\, {s}_{7,1} + {s}_{7,2} + 2\, {s}_{7,3}$ & ${s}_{9,1} + 2\,
{s}_{9,2}$ & $0$ &  \\ \hline
${s}_{5,3}$ & $2\, {s}_{7,1}$ & ${s}_{9,1} + 2\, {s}_{9,2}$ & $0$ &  \\
\hline
${s}_{5,4}$ & $2 {s}_{7,1} + 2 {s}_{7,2} + {s}_{7,3} + 2 {s}_{7,4}$ & $0$ & $%
0$ &  \\ \hline
${s}_{6,2}$ & $2\, {s}_{8,1} + 2\, {s}_{8,2} + 2\, {s}_{8,3}$ & $0$ & $2\, {s%
}_{12,1}$ &  \\ \hline
${s}_{6,3}$ & $2\, {s}_{8,2} + {s}_{8,3}$ & ${s}_{10,1}$ & $0$ &  \\ \hline
${s}_{6,4}$ & $2\, {s}_{8,1} + {s}_{8,2} + {s}_{8,4}$ & $0$ & $2\, {s}%
_{12,1} $ &  \\ \hline
${s}_{6,5}$ & $2\, {s}_{8,1} + {s}_{8,3}$ & $2\, {s}_{10,1}$ & $2\, {s}%
_{12,1}$ &  \\ \hline
${s}_{7,1}$ & ${s}_{9,1} + 2\, {s}_{9,2}$ & $0$ &  &  \\ \hline
${s}_{7,2}$ & $2\, {s}_{9,2} + 2\, {s}_{9,3}$ & $2\, {s}_{11,1}$ &  &  \\
\hline
${s}_{7,3}$ & $2\, {s}_{9,2} + 2\, {s}_{9,3}$ & $2\, {s}_{11,1}$ &  &  \\
\hline
${s}_{7,4}$ & $2\, {s}_{9,1} + {s}_{9,2}$ & $0$ &  &  \\ \hline
${s}_{8,1}$ & $2\, {s}_{10,1} + 2\, {s}_{10,2}$ & ${s}_{12,1}$ &  &  \\
\hline
${s}_{8,2}$ & ${s}_{10,1} + 2\, {s}_{10,2}$ & ${s}_{12,1}$ &  &  \\ \hline
${s}_{8,3}$ & $2\, {s}_{10,2}$ & ${s}_{12,1}$ &  &  \\ \hline
${s}_{8,4}$ & ${s}_{10,1}$ &  &  &  \\ \hline
${s}_{9,1}$ & $2\, {s}_{11,1}$ &  &  &  \\ \hline
${s}_{9,2}$ & $2\, {s}_{11,1}$ &  &  &  \\ \hline
${s}_{10,2}$ & ${s}_{12,1}$ &  &  &  \\ \hline
\end{tabular}
}

\begin{tabular}{|l|l|l|}
\multicolumn{3}{c}{\textbf{Table B${}_2$.} \textsl{\ Nontrivial }$%
P^{k}(s_{r,i})$ \textsl{for $p=5$}} \\ \hline
$s_{r,i}$ & $\mathcal{P}^{1}(s_{r,i})$ & $\mathcal{P}^{2}(s_{r,i})$ \\ \hline
${s}_{1,1}$ & ${s}_{5,1} + 4\, {s}_{5,3}$ & $0$ \\ \hline
${s}_{2,1}$ & $4\, {s}_{6,1} + {s}_{6,2} + 4\, {s}_{6,4}$ & ${s}_{10,1}$ \\
\hline
${s}_{2,2}$ & ${s}_{6,1} + {s}_{6,2} + 4\, {s}_{6,4}$ & ${s}_{10,1}$ \\
\hline
${s}_{3,1}$ & $4\, {s}_{7,1} + {s}_{7,3}$ & ${s}_{11,1}$ \\ \hline
${s}_{3,2}$ & ${s}_{7,2} + {s}_{7,3} + 4\, {s}_{7,4}$ & $2\, {s}_{11,1}$ \\
\hline
${s}_{3,3}$ & ${s}_{7,1} + {s}_{7,2} + 4\, {s}_{7,4}$ & ${s}_{11,1}$ \\
\hline
${s}_{4,1}$ & $4\, {s}_{8,1} + {s}_{8,2} + {s}_{8,3}$ & ${s}_{12,1}$ \\
\hline
${s}_{4,2}$ & ${s}_{8,2} + {s}_{8,3} + 4\, {s}_{8,4}$ & $2\, {s}_{12,1}$ \\
\hline
${s}_{4,3}$ & ${s}_{8,1} + 4\, {s}_{8,2} + {s}_{8,3}$ & ${s}_{12,1}$ \\
\hline
${s}_{4,4}$ & ${s}_{8,1} + 4\, {s}_{8,4}$ & ${s}_{12,1}$ \\ \hline
${s}_{5,1}$ & $4\, {s}_{9,1} + {s}_{9,2}$ &  \\ \hline
${s}_{5,2}$ & ${s}_{9,2} + {s}_{9,3}$ &  \\ \hline
${s}_{5,3}$ & $4\, {s}_{9,1} + {s}_{9,2}$ &  \\ \hline
${s}_{5,4}$ & ${s}_{9,1} + {s}_{9,3}$ &  \\ \hline
${s}_{6,2}$ & ${s}_{10,1} + {s}_{10,2}$ &  \\ \hline
${s}_{6,3}$ & $2\, {s}_{10,2}$ &  \\ \hline
${s}_{6,4}$ & $4\, {s}_{10,1} + {s}_{10,2}$ &  \\ \hline
${s}_{6,5}$ & ${s}_{10,1} + 4\, {s}_{10,2}$ &  \\ \hline
${s}_{7,2}$ & ${s}_{11,1}$ &  \\ \hline
${s}_{7,3}$ & $2\, {s}_{11,1}$ &  \\ \hline
${s}_{7,4}$ & $4\, {s}_{11,1}$ &  \\ \hline
${s}_{8,3}$ & $2\, {s}_{12,1}$ &  \\ \hline
${s}_{8,4}$ & $3\, {s}_{12,1}$ &  \\ \hline
\end{tabular}

\bigskip \bigskip

\textbf{References}
\end{center}

[BGG] I. N. Bernstein, I. M. Gel'fand and S. I. Gel'fand, Schubert cells and
cohomology of the spaces G/P, Russian Math. Surveys 28 (1973), 1-26.

[BH] A. Borel and F. Hirzebruch, Characteristic classes and homogeneous
spaces (I), Amer. J. Math. 80 (1958), 458--538.

[BS$_{1}$] R. Bott and H. Samelson, The cohomology ring of G/T, Proc. Nat.
Acad. Sci. U. S. A. 41 (1955), 490--493.

[BS$_{2}$] R. Bott and H. Samelson, Application of the theory of Morse to
symmetric spaces, Amer. J. Math., Vol. LXXX, no. 4 (1958), 964-1029.

[BSe$_{1}$] A. Borel and J. P. Serre, D\'{e}termination des $p$-puissances r%
\'{e}duites de Steenrod dans la cohomologie des groupes classiques.
Applications, C. R. Acad. Sci. Paris 233, (1951). 680--682.

[BSe$_{2}$] A. Borel and J. P. Serre, Groupes de Lie et puissances r\'{e}%
duites de Steenrod, Amer. J. Math. 75(1953), 409-448.

[Ch] C. Chevalley, Sur les D\'{e}compositions Celluaires des Espaces G/B, in
Algebraic groups and their generalizations: Classical methods, W. Haboush
ed. Proc. Symp. in Pure Math. 56 (part 1) (1994), 1-26.

[D] J. Dieudonn\'{e}, A history of Algebraic and Differential Topology,
1900-1960, Boston; Basel, 1989.

[D$_{1}$] H. Duan, The degree of a Schubert variety, Adv. in Math.,
180(2003), 112-133.

[D$_{2}$] H. Duan, Self-maps of the Grassmannian of complex structures,
Compositio Math., 132 (2002), 159-175.

[D$_{3}$] H. Duan, Multiplicative rule of Schubert classes,
Invent. Math., 159 (2)(2005), 407-436.

[D$_{4}$] H. Duan, On the inverse Kostka matrix, J. Combinatorial Theory A,
103(2003), 363-376.

[DZ] H. Duan and Xuezhi Zhao, Algorithm for multiplying Schubert classes,
arXiv: math.AG/0309158.

[DZZ] H. Duan, Xu-an Zhao and Xuezhi Zhao, The Cartan matrix and enumerative
calculus, J. Symbolic computation, 38(2004), 1119-1144.

[E] C. Ehresmann, Sur la topologie de certains espaces homogenes, Ann. of
Math. (2) 35 (1934), 396--443.

[GH] P. Griffith and J. Harris, Principles of algebraic geometry, Wiley, New
York, 1978.

[H] H.C. Hansen, On cycles in flag manifolds, Math. Scand. 33 (1973),
269-274.

[Hu] J. E. Humphreys, Introduction to Lie algebras and representation
theory, Graduated Texts in Math. 9, Springer-Verlag New York, 1972.

[K] S. Kleiman, Problem 15. Rigorous fundation of the Schubert's enumerative
calculus, Proceedings of Symposia in Pure Math., 28 (1976), 445-482.

[La] T. Lance, Steenrod and Dyer-Lashof operations on BU, Trans. Amer. Math.
Soc. 276(1983), 497-510.

[Le] C. Lenart, The combinatorial of Steenrod operations on the cohomology
of Grassmannians, Adv. in Math. 136(1998), 251-283.

[M] I. G. Macdonald, Symmetric functions and Hall polynomials, Oxford
Mathematical Monographs, Oxford University Press, Oxford, second ed., 1995.

[P] F. Peterson, A mod--p Wu formula, Bol.Soc.Mat.Mexicana 20(1975), 56-58.

[S] B. Shay, mod--p Wu formulas for the Steenrod algebra and the Dyer-Lashof
algebra, Proc. AMS, 63(1977), 339-347.

[So] F. Sottile, Four entries for Kluwer encyclopaedia of Mathematics,
arXiv: Math. AG/0102047.

[SE] N. E. Steenrod and D. B. A. Epstein, Cohomology Operations, Ann. of
Math. Stud., Princeton Univ. Press, Princeton, NJ, 1962.

[St] N. E. Steenrod, Cohomology operations, and obstructions to extending
continuous functions, Adv. in Math. 8 (1972), 371--416.

[Su] T. Sugawara, Wu formulas for the mod--3 reduced power operations, Mem.
Fac. Sci. Kyushu Univ. Ser. A 33(1979), 297-309.

[SW] N.E. Steenrod and J.H.C. Whitehead, Vector fields on the $n$--sphere,
Proc. Nat. Acad. Sci. U. S. A. 37, (1951). 58--63.

[W] R. Wood, Problems in the Steenrod algebra, Bull. London Math. Soc.
30(1998), 449-517.

[Wu] T. Wu, Les $i$--carr\'{e}s dans une vari\'{e}t\'{e} grassmanniene, C.
R. Acad. Sci. Paris 230 (1950), 918--920.

\end{document}